\documentclass[10pt, a4paper]{article}

\usepackage[cp850]{inputenc}

\usepackage[T1]{fontenc}

\usepackage{graphicx,color}

\def\reel{\hbox{I\hskip -2pt R}}

\def\nne{\hbox{I\hskip -2pt N }}
\def\esp{\hbox{I\hskip -2pt E }}
\def\a{{\bf a}}

\def\b{{\bf b}}

\def\<{\langle}

\def\>{\rangle}

\newtheorem{cor}{Corollary}

\newtheorem{prop}{Proposition}

\title{A Random Difference Equation with Dufresne Variables revisited}

\author{Jean-Fran\c{c}ois Chamayou \thanks{Laboratoire de Statistique et Probabilit\'es, Universit\'e Paul Sabatier, 31062 Toulouse, France. e-mail: \texttt{jfchamay@math.univ-toulouse.fr}}}

\date{\today}

\begin{document}

\maketitle

\begin{abstract}
The Dufresne laws (laws of product of independent random variables with gamma and beta distributions) occur as stationary distribution of certain Markov chains $ X_n $ on $\reel$  defined by:
\begin{equation}
 X_n = A_n ( X_{n-1}  + B_n ) 
 \end{equation} 
 where $ X_0 , (A_1,B_1),...,(A_n,B_n) $ are independent and the $(A_i,B_i)'$s are identically distributed. 
 This paper generalizes an explicit example where $A$ is the product of two independent $\beta_{a,1} , \beta_{b,1} $ and  $B \sim \gamma_1 $ or $ \gamma_2 $.\\



\end{abstract}

\section{Introduction.}
Let us recall that the real  Dufresne laws are distributions on $(0,\infty)$ defined as follows: 
Let $p$ and $q$ be in $ \nne = {0,1,2,...} $ and $ \a = ( a_1,...,a_p)$ and $ \b =(b_1,...,b_q)$ be two sequences of positive numbers. We write for $s \geq 0$ and $p > 0 $
$$ ( \a)_s = \prod_{j=1}^{p}  \frac{\Gamma (a_j + s)}{\Gamma(a_j)} $$
and $ (\a)_s = 1 $ for $ p=0 $. As introduced  in [3,6], the Dufresne law 
$ D(\bf{a};\bf{b}) $ on $ (0, + \infty) $ (if it exists) is defined by its Mellin transform:
$$ \int_{0}^{\infty} x^s D(\a,\b) (dx) = \frac{(\a)_s}{(\b)_s} $$
where $s$ is such that $a_j+s$ and $b_i+s$ are positive for all $i$ and $j.$ 
 For $ p = q = 1 $ an example is given the beta distribution defined for $ a > 0 $ and $ b > 0 $ by
\begin{equation}
\beta_{a,b}(dx)=\frac{1}{B(a,b)}x^{a-1}(1-x)^{b-1}{\bf 1}_{(0,1)}(x)dx,
\end{equation}
where $B(a,b)= \frac{\Gamma(a)\Gamma(b)}{\Gamma(a+b)} .$
and which therefore can also be written as $ D(a;a+b). $ 
Similarly for $ p=1 $ and $ q=0 $ 
the gamma distribution with scale parameter $1$ and  shape parameter $c>0$ is the probability on $(0,\infty)$ defined by

\begin{equation}
\gamma_c(dx)=\frac{1}{\Gamma(c)} x^{c-1} e^{-x}{\bf 1}_{(0,\infty)}(x)dx
\end{equation}
is the Dufresne distribution $ D(c,-) $.
\section{General case}
Let us now develop an extension to paragraph 5 in [3]:
\begin{prop}
Let  $A,B$ and $X$ be independent positive random variables such that $A$ has distribution $D(\alpha,\beta; \alpha+1,\beta+1)$ (i.e. is the distribution of product of two independent random variables $ \beta_{\alpha,1}, \beta_{\beta,1}$) and such that $ B $ has distribution $ \gamma_u = D(u;-) $, where $u$ is a positive real parameter. Then $X$ and $ A ( X + B) $ have the same distribution if and only if $ \Phi (s) = \esp{(e^{-sX} )}  $ is the unique continuous solution for $ s \geq 0 $ with $ \Phi(0) = 1 $ of the following differential equation:
\begin{equation}
 (s+1)^u s \Phi''(s) + (\alpha + \beta +1) (s+1)^u \Phi'(s) + \frac{\alpha \beta}{s} ((s+1)^u - 1) \Phi(s) = 0 
 \end{equation}
\end{prop}

\noindent\textbf{Proof}. Suppose that  $X$ and $ A(X+B) $ have the same distribution and consider  $ \Phi (s) = \esp{(e^{-sX} )}  $
for $ s \geq 0 $. The law of $A$ is:
$$ \frac{ \alpha \beta }{\beta - \alpha } [ t^{\alpha-1} - t^{\beta-1} ]\ {\bf 1}_{(0,1)}(t) dt $$
If $B$ has the gamma distribution $ D(u;-) $, we have
$$ \esp{(e^{-sB}) } = \frac{1}{(1+s)^u} $$
and the functional equation can be written:
\begin{equation}
\Phi(s)=\esp{(e^{-sA(X+B)})} =\esp{(\esp{(e^{-sA(X+B)}\mid A)})} =\esp{(\frac{\Phi(sA)}{(1+sA)^u})} =
\alpha \beta \int_{0}^{1} \int_{0}^{1} \frac{\Phi(sab) a^{\alpha-1}b^{\beta-1}}{(1+sab)^u} da db
\end{equation}
leading to: 
\begin{equation}
\Phi(s)=\frac{\alpha \beta}{s^\alpha}\int_{0}^{s} \int_{0}^{1} \frac{\Phi(tb) t^{\alpha-1} b^{\beta-1}}{(1+tb)^u} dt db
\end{equation}

we rewrite this as
\begin{equation}
s^\alpha \Phi(s) = \alpha \beta \int_{0}^{s} \int_{0}^{t} \frac{\Phi(v) t^{\alpha-1} v^{\beta-1}}{t^{\beta} (1+v)^u} dt dv
\end{equation}
Taking the derivative of both sides gives:
\begin{equation}
(s^\alpha \Phi(s))' = \alpha \beta s^{\alpha -\beta -1} \int_{0}^{s} \frac{\Phi(v) v^{\beta-1}}{(1+v)^u} dv
\end{equation}
and this easily leads to the differential equation (2.4). To show the uniqueness of $\Phi$ , an elementary way is to see that all other solution of (2.4) in a neighborhood of zero can be written as $ s \longmapsto \Phi(s) z(s) $, where the unknown function $z$ satisfies 
$$ s \Phi(s) z''(s) + (2 \Phi'(s) +( \alpha + \beta +1) \Phi(s)) z'(s) = 0 $$
Since from (2.4) $ \Phi'(0) =- \frac{u \alpha \beta}{\alpha + \beta +1}$, it is easily seen that $z$ is continuous on 0 with $z(0)=1$ iff $z$ is a constant.

\section{Case $u=1$}
Let us study first the case where $u=1$. Then the problem becomes a particular case of the Dufresne result $\alpha, \beta , c=1 $ (see Proposition 3 in [3] see also [5]) and the distribution of $X$ is $ D(\alpha,\beta; \alpha + \beta +1) $ . A direct proof of this result can also been obtained from the differential equation (2.4), which becomes for $u=1$ :
$$ (s+1) s \Phi''(s) + (\alpha + \beta +1) (s+1)  \Phi'(s) + \alpha \beta \Phi(s) = 0 $$  
whose solutions  can be calculated by the standard tools of hypergeometric differential equations.
As noted in [3] the explicit expression of the Dufresne distribution in terms of Whittaker confluent hypergeometric function $ W(p,q;x) $ is:
\begin{equation}
\frac{\Gamma( \alpha + \beta +1)}{\Gamma(\alpha)\Gamma(\beta)} e^{-x/2} x^{\frac{\alpha+ \beta-3}{2}} W(-\frac{1+\alpha+\beta}{2},\frac{\beta-\alpha}{2};x){\bf 1}_{(0,\infty)}(x)dx 
\end{equation}
 (see Erdelyi et al. 5.19 [8], or Gradsteyn and Ryzhik 7.621(6) [11]).
 
\section{Case $u=2$}
We now consider the more difficult case $ u=2 $.  Proposition 12 in [3] was  curiously  involving the golden ratio $ \frac{1+ \sqrt{5}}{2} $. We extend it as follows:

\begin{prop}
Let $A,B$ and $X$ be independent positive random variables such that A has distribution $D(\alpha,\beta; \alpha+1,\beta+1)$ and such that $B$ has distribution  $ \gamma_2 = D(2;-) $. Then $X$ and $A(X+B)$ have the same distribution if and only if $X$ has the distribution:
 $$ D(\alpha +\rho,\beta +\rho ;  \alpha + \beta + 1) * D(-\rho , -) $$
 where $ \rho = \frac{1 - \sqrt{1+4\alpha \beta}}{2} $ and the operation $*$ stands for an additive convolution of Dufresne variables.
 
 \end{prop} 
\noindent\textbf{Proof}. For $n=2$ the differential equation (2.4) becomes:

\begin{equation}
 (s+1)^2 s \Phi''(s) + (\alpha+\beta+1) (s+1)^2 \Phi'(s) + \alpha \beta (s+2) \Phi(s) = 0 
 \end{equation} 
The characteristic equation of this differential equation relative to the singular point $s=-1$ is $ p^2 -p -\alpha \beta $, whose solutions are $ \rho_+ =  \frac{1 + \sqrt{1+4\alpha \beta}}{2} $ or $ \rho =  \frac{1 - \sqrt{1+4\alpha \beta}}{2} $. Changing the function $ \Phi $ into the new unknown function $ z$ defined by $ \Phi (s) = (1+s)^\rho z(s) $ leads to the following hypergeometric differential equation for $z$:
$$ s(1+s) z''(s) + ( \alpha + \beta +1 + ( 2 \rho + \alpha + \beta + 1 )s )z'(s) + (\rho + \alpha )(\rho + \beta ) z(s) = 0 $$
whose unique solution continuous on $ [ 0 , \infty ) $ and equal to $1$ on $0$ is the Gauss  hypergeometric function :
$$   _2F_1 ( \rho + \alpha ,\rho + \beta ;  \alpha + \beta + 1 ; -s ) $$
(see Rainville p. 54 [12]). Thus
$$ \Phi(s) = (1+s)^\rho {}_2F_1 ( \rho + \alpha ,\rho + \beta ; \alpha + \beta + 1 ; -s ) $$
We know that the Laplace transform of the Dufresne distribution $D ( \rho + \alpha ,\rho + \beta ;  \alpha + \beta + 1 )$ is $  _2F_1 ( \rho + \alpha ,\rho + \beta ; \alpha + \beta + 1 ; -s ) $ and that the Laplace transform of the distribution $ \gamma_{-\rho} $ is $ (1+s)^\rho $: the result is proved.
Note that the positivity of the parameters is given by the condition 
$ min(\alpha, \beta) +1 > max(\alpha , \beta)$, this condition shows the exchange of the product of two beta variates: $D(\alpha,\beta; \alpha+1,\beta+1)$ is equivalent to $D(\alpha,\beta; \beta+1,\alpha+1)$.

\section{Corollaries}
As a consequence, we get a possibly new identity for the generalized hypergeometric function $ _3F_2 $.
\begin{cor}
For all integers $ n>0 $ we have the following identity:
\begin{equation}
(-\rho)_n {}\ _3F_2 (\rho + \alpha ,\rho + \beta,-n;  \alpha + \beta  +1 , \rho -n ; 1 ) = 
\frac{1}{(n+1)^2} (2- \rho )_n  {}\  _3F_2 (\rho + \alpha ,\rho + \beta,-n; \alpha+\beta + 1 , \rho -3 -n ; 1)
\end{equation}
\end{cor}
\noindent\textbf{Proof}. With the notations of the previous proposition, we use the fact that $ \esp{(X^n)} = \esp{(A^n)} \esp{((X+B)^n)}$. Now the distribution of $ X+B $ is $ D(\alpha +\rho,\beta +\rho ; \alpha+ \beta + 1) * D(-\rho + 2 , -) $, from proposition (2) and the convolution properties of the gamma distributions. We can now apply  Proposition 7 in [3] for computing the moments of $ X$ and $ X+B $. Since $ \esp{(A^n)} = \frac{1}{(n+1)^2} $
, the corollary is proved.\\
Another consequence, where we get a possibly new integral formula for Gauss hypergeometric functions, arises from the application of the Laplace transform to the result.
\begin{cor}
If $ \alpha , \beta $ are  positive, if  $\rho=  \frac{1 - \sqrt{1+4\alpha \beta}}{2} $ and if $ _2F_1( a,b;c;s) $ is  the Gauss hypergeometric function then we have :
\begin{equation}
(1+s)^\rho {}_2F_1 ( \rho + \alpha ,\rho + \beta ;  \alpha + \beta + 1 ; -s )=
\int_{0}^{1} \int_{0}^{1}\frac{\alpha \beta x^{\alpha-1} y^{\beta-1}}{(1+sxy)^{-\rho+2}} {}_2F_1 ( \rho + \alpha ,\rho + \beta ;\alpha+\beta + 1 ; -sxy ) dxdy
 \end{equation}
 \end{cor}
\noindent\textbf{Proof}. With the notations of the previous proposition, we use the fact that $\esp{(e^{-sX})}= \esp{(e^{-sA(X+B)})} $. Now the Laplace transform of the distribution of $X+B$ is $ (1+s)^{\rho-2} {}_2F_1 ( \rho + \alpha ,\rho + \beta ;\alpha +\beta + 1 ; -s ). $ From Proposition 2, from the convolution properties of the gamma function and since $A$ is the product of two independent $ \beta_{\alpha,1}, \beta_{\beta,1}$ the corollary is proved.                                        

\section{Application to the triggered shot noise}
In physics and in insurance mathematics (see [5])the following process occurs:
\begin{equation}
  Z(t) = \sum_{0 \leq T_i \leq t} B_i e^{-p(t-T_i)}
  \end{equation} 
where the $(T_i)_{i \in \nne^* }$ are the arrival dates of a Poisson process with rate $ \lambda $ and the $(B_i)_{i \in \nne^*}$  are i.i.d. random variables . $Z(t)$ represents, for example, the shot noise as in ref. [15]. The stationary limit of this process exists if $ E(\log |B_i|) < \infty $ and is the same as the distribution limit of the sequence of random variables $(X_i)_{i \in \nne^*} $ defined by 
$$ X_1 = A_1 B_1 , X_{m+1} = A_{m+1} ( X_m + B_{m+1} ) , m> 0 $$
the $ A_j $ being powers in $ \frac{p}{\lambda}$ of independent random variables uniform on $[0,1]$ .\\
Now consider a triggered shot noise with period k (see [1,2,4]):
$$ Z_k(t) = \sum_{i| 0 \leq T_{ki} \leq t} B_i e^{-(\sum_{j=2}^{k} p_{(j-1)}(T_{(k(i-1)+j)}-T_{(k(i-1)+(j-1)})+p_{(k)}(t-T_{(ki)}))} $$
The k sucessive Poisson events are  perturbating  the exponential decay differently. The same arguments as in [9] show that the distribution limit of $ Z_{T_{km}} $ is given by the limit of $X_m$ , but the $A_i$ are independent with the same distribution as $\Pi_{j=1}^{k} U_j ^{p_{j}/\lambda} $ the $U_j$ being independent random variables uniform on $[0,1]$
which correspond to the studied cases $ u=1,2 $ where $ \alpha=  \frac{p_1}{\lambda}, \beta=  \frac{p_2}{\lambda} $  and the $ B_i$ are $ \gamma_1 $ or $ \gamma_2 $ distributed .

 \section{Products of more than 2 betas}
\begin{prop}
Let  $A,B$ and $X$ be independent positive random variables such that $A$ has distribution $D(\alpha,\beta,\gamma; \alpha+1,\beta+1,\gamma+1)$ (i.e. is the distribution of product of three independent random variables $ \beta_{\alpha,1}, \beta_{\beta,1}, \beta_{\gamma,1}$) and such that $ B $ has an exponential distribution $ \gamma_1 = D(1;-) $. Then $X$ and $ A ( X + B) $ have the same distribution if and only if $ \Phi (s) = \esp{(e^{-sX} )}  $ is the unique continuous solution for $ s \geq 0 $ with $ \Phi(0) = 1 $ of the following differential equation:
\begin{equation}
 (1+s) s^2 \Phi'''(s) + (\alpha + \beta +\gamma+3) (1+s) s \Phi''(s) +(1+s)(\alpha\beta+\alpha\gamma+\beta\gamma+\alpha+\beta+\gamma+1) \Phi'(s)+ \alpha \beta \gamma \Phi(s) = 0 
 \end{equation}
\end{prop}
\noindent\textbf{Proof}.
Using the same arguments as in proposition (2), we get
\begin{equation}
\Phi(s)=\alpha \beta\gamma \int_{0}^{1} \int_{0}^{1}\int_{0}^{1}  \frac{\Phi(sabc) a^{\alpha-1}b^{\beta-1}c^{\gamma-1}}{(1+sabc)} da db dc
 \end{equation}
from which after 3 differentiations we derive the differential equation. The problem  again becomes  a particular case of Dufresne distribution   $ D(\alpha,\beta,\gamma; c_1,c_2) $ since the differential equation is of generalized hypergeometric type
with solution $$  {}_3F_2 (  \alpha ,\beta ,\gamma; c_1,c_2 ; -s ) $$ see [14] the parameters $c_1, c_2$ being the roots of the second order equation generated by the non linear system given by the identification of the parameters of equation (7.14) with those of the canonical equation in [14]. i.e.:
$$ c_1 + c_2 +1= \alpha+\beta+\gamma+3$$
and
$$ c_1 c_2 = \alpha\beta+\alpha\gamma+\beta\gamma+\alpha+\beta+\gamma+1 $$
The roots are complex conjugate (or real positive),
$$ c_1 = \frac{\alpha+\beta+\gamma+2 +  i \sqrt{4(\alpha\beta+\alpha\gamma+\beta\gamma+\alpha+\beta+\gamma+1)-(\alpha+\beta+\gamma+2)^2}}{2}$$
 the case $\alpha=\beta=\gamma$ has been studied by Dufresne in [7] where he justifies the fact of having complex parameters for the product of beta distributions, such complex parameters have appeared already for the distribution of triggered shot noise [2], it is easily verified that the real part and the module of the roots (or the  two real roots ) are greater than the elements of $\overline{(\alpha ,\beta ,\gamma)\bigcap  max(\alpha ,\beta ,\gamma)}$  to build the Laplace transform of a beta variable. \\
Products of $n >3$  betas could be handled in the same way, leading again to Dufresne variables but with increasing complexity to calculate the parameters since for $n>5$ it will be necessary to use iterative numerical methods to solve the system of $n-1$
non-linear equations to evaluate the $n-1$ unknown parameters.
 For instance for $n=4$ the differential equation is:
\begin{eqnarray*}
(1+s) s^3 \Phi^{(iv)}(s) +(\alpha+\beta +\gamma+\delta+6) (1+s) s^2 \Phi'''(s) \\ +(\alpha\beta+\alpha\gamma+\alpha\delta+\beta\gamma+\beta\delta+\gamma\delta+
3(\alpha+\beta+\gamma+\delta)+7) (1+s)s \Phi''(s) \\
+(\alpha\beta\gamma+\alpha\gamma\delta+\alpha\beta\delta+\beta\gamma\delta+
\alpha\beta+\alpha\gamma+\alpha\delta+\beta\gamma+\beta\delta+
\gamma\delta+\alpha+\beta+\gamma+\delta+1) (1+s) \Phi'(s)\\
+\alpha \beta \gamma\delta \Phi(s) = 0 
\end{eqnarray*}
the solution is:
$$ {}_4F_3 (  \alpha ,\beta ,\gamma, \delta; c_1,c_2, c_3 ; -s ) $$
the unkown parameters are the roots of the third order equation resulting from the 3 non linear equations due to the identification with the parameters of the fourth order hypergeometric equation given in [13]i.e.: 
$$ c_1 + c_2 + c_3 +3 = \alpha+\beta+\gamma+\delta +6 $$
$$c_1 + c_2 + c_3 + c_1 c_2 + c_2 c_3 + c_3 c_1+1 =\alpha\beta+\alpha\gamma+\alpha\delta+\beta\gamma+\beta\delta+\gamma\delta+
3(\alpha+\beta+\gamma+\delta)+7$$
$$ c_1 c_2 c_3 = \alpha\beta\gamma+\alpha\gamma\delta+\alpha\beta\delta+\beta\gamma\delta+
\alpha\beta+\alpha\gamma+\alpha\delta+\beta\gamma+\beta\delta+
\gamma\delta+\alpha+\beta+\gamma+\delta+1 $$
There are 3 real positive roots or one real positive root and two complex conjugate roots which is allowed according to the Dufresne results [7] the positivity being given by the Descartes/Laguerre rule of signs, see appendix.\\
For the case $n=5$ the solution is:
$$ {}_5F_4 (  \alpha ,\beta ,\gamma, \delta, \eta; c_1,c_2, c_3, c_4 ; -s ) $$
the unkown parameters are the roots of the fourth order equation resulting from the 4 non linear equations due to the identification with the parameters of the fith order hypergeometric equation :
$$ c_1 + c_2 + c_3 + c_4 +6 = \alpha+\beta+\gamma+\delta+ \eta + 10 $$
\begin{eqnarray*}
3(c_1 + c_2 + c_3 + c_4 )+ c_1 c_2 + c_2 c_3 + c_3 c_1 + c_1 c_4 + c_2 c_4+ c_3 c_4+7=\\ \alpha\beta+\alpha\gamma+\alpha\delta+\beta\gamma+\beta\delta+\gamma\delta
+\alpha\eta+\beta\eta+\gamma\eta+\delta\eta
+6(\alpha+\beta+\gamma+\delta+\eta)+25
\end{eqnarray*}
\begin{eqnarray*}
c_1c_2c_3+c_1c_2c_4+c_1c_3c_4+c_2c_3c_4+c_1+c_2+c_3+c_4+c_1c_2+c_2c_3+c_3c_1+c_1c_4+c_2c_4+c_3c_4+1=\\ \alpha\beta\gamma+\alpha\gamma\delta+\alpha\beta\delta+\beta\gamma\delta
+\alpha\beta\eta+\alpha\gamma\eta+\alpha\delta\eta+\beta\gamma\eta+\beta\delta\eta+\gamma\delta\eta\\
+3(\alpha\beta+\alpha\gamma+\alpha\delta+\beta\gamma+\beta\delta+\alpha\eta+\beta\eta+\gamma\eta
+\gamma\delta+\delta\eta)+7(\alpha+\beta+\gamma+\delta+\eta)+15
\end{eqnarray*}
\begin{eqnarray*} 
c_1 c_2 c_3 c_4 =
\alpha\beta\gamma\eta+\alpha\gamma\delta\eta+\alpha\beta\delta\eta+\beta\gamma\delta\eta
+\alpha\beta\gamma\delta\\
\alpha\beta\gamma+\alpha\gamma\delta+\alpha\beta\delta+\beta\gamma\delta
+\alpha\beta\eta+\alpha\gamma\eta+\alpha\delta\eta+\beta\gamma\eta+\beta\delta\eta+\gamma\delta\eta\\
+(\alpha\beta+\alpha\gamma+\alpha\delta+\beta\gamma+\beta\delta+\alpha\eta+\beta\eta+\gamma\eta
+\gamma\delta+\delta\eta)+(\alpha+\beta+\gamma+\delta+\eta)+1
\end{eqnarray*}
There are 4 real positive roots or two real positive roots and two complex conjugate roots or two pairs of complex conjugate roots which is allowed according to the Dufresne results [7],the positivity being given by the Descartes/Laguerre rule of signs, see appendix.
The process could be extended by induction using an automatic procedure, the main difficulty being then for $n> 5$ to get a good initial guess of the location of the solution of the non-linear system by an iterative method, for the sake of simplicity the choice of the set of $\overline{(\alpha ,\beta ,\gamma, \delta, \eta,...)\bigcap  max(\alpha ,\beta ,\gamma, \delta, \eta , ...)}$ components for the starting array of the real roots of the iteration process could be made since their positivity is assumed from the Descartes/Laguerre rule of signs and each real positive  parameters $ c_i $ must be greater than one of the associated $ \alpha ,\beta ,\gamma, \delta, \eta , ...$ to build the Laplace transform of a beta variable. Otherwise the calculations could be made sequentially: set the parameters $(\alpha ,\beta ,\gamma, \delta, \eta,...)$ in increasing order $ \alpha_i, i= 1 , ... ,n+1$ 
compute the n roots starting the iterations  with the initial guess $ \alpha_1 $ and then $  max( c_{i-1} , \alpha_i) , i=2,...,n $. For instance, if
$\alpha=\beta=\gamma=\delta=\eta=\zeta=1$ then 
 $$c_1 = 3,c_2 = \frac{3 + i\sqrt{3}}{2} ,c_3 = \frac{3 - i\sqrt{3}}{2},c_4 = \frac{5 + i\sqrt{3}}{2},c_5 = \frac{5 - i\sqrt{3}}{2}  $$ 
 (using wolframalpha) equivalent results have been tested which show that the special case: equality to 1 of the all the parameters gives a set of real integer coefficients for a quasi cyclotomic polynomial, the roots are the roots of unit (the real root 1 is lacking ) the unit circle being centered at $(2,0)$ on the complex plane ( C++ program available on request, see also Sloane [14]). The corresponding cycltomic polynomial is the following according to the identity [13]: 
 $$ \sum_{m=0}^{m=n} (-1)^m y^{n-m} $$
 An extension to the case $\alpha=\beta=\gamma=\delta=\eta=\zeta=...$ gives the roots on the circle of radius $\alpha $ centered at $ (\alpha+1,0) $ in the complex plane. This extension due to the generalization of the relation of theorem 1.1 in [13]
 $$ \sum_{j=r+1}^{n} C_{n}^{j} C_{j-1}^{r} \alpha^{j-1-r}  =  \sum_{j=r+1}^{n} C_{j-1}^{r} (1+ \alpha)^{j-1-r} ,  0\leq r \leq n-1 $$
 For instance if $\alpha = 2$ the polynomial is $ x^3-11 x^2+43 x-65 $ the roots are $ x = 3-2 i, x = 3+2 i,  x = 5$ 
on the circle of radius 2 centered at $(3,0)$ the real root $ 1 $ is lacking see figure.If all the parameters are equal to $ \frac{1}{n} $ at the limit $n$ tending to infinity it can be conjectured that the result is degenereted to a circle point of radius $ \frac{1}{n} $ centered at $ (1+\frac{1}{n} , 0)$ showing an infinite product of $ \beta_{ \frac{1}{n},1} $.

\newpage

\section{References}\begin{enumerate}

\item  \textsc{\sc Blanc-Lapierre M.A.} (1955) `Consid\'{e}ration sur certains processus ponctuels et sur les fonctions al\'{e}atoires associ\'{e}es ' {\em Colloque Analyse Statistique, Centre Belge de Recherches Math\'ematiques} 25-55.

\item  \textsc{\sc Chamayou, J.-F.} (1996) `A Case of a Random Difference Equation Connected with a Filtered Binary Process' {\em Adv. Appl. Math.} {\bf 17} 88-100.

\item  \textsc{\sc Chamayou, J.-F., Letac, G.} (1999) `Additive Properties of the Dufresne Laws and Their
Multivariate Extension' {\em J. Theor. Prob.} {\bf 12} 1045-1066.

\item  \textsc{\sc Chamayou, J.-F., Dunau J.L.} (2002) `Random difference equations with logarithmic distribution and the triggered shot noise' {\em Adv. Appl. Math.} {\bf 29} 454-470.

\item  \textsc{\sc Dufresne D.} (1990) 'The distribution of a perpetuity, with applications to risk theory and pension funding' {\em Scand. Actuarial J.} {\bf 1/2} 39-79.

\item  \textsc{\sc Dufresne D.} (1996) `On the stochastic equation $ \it{L} (X) = \it{L} (B(X+C)) $ and a property of gamma distributions' {\em Bernoulli} {\bf 2} 287-291.
 
\item  \textsc{\sc Dufresne D.} (2010) `The Beta Product Distribution with Complex Parameters' {\em Communications Statistics Theory Methods} {\bf 39} 837-854.
\item  \textsc{\sc Erdelyi A., \it{et al.}} (1953a) {\em Higher Transcendental Functions}, McGraw Hill, New York.

\item  \textsc{\sc Erdelyi A., \it{et al.}} (1953b) {\em Tables of Integral Transforms}, McGraw Hill, New York.

\item  \textsc{\sc Feller W.} (1966) {\em An Introduction to Probability Theory and its Applications, Vol. II}, Wiley, New York.

\item  \textsc{\sc Gradstein, I.S., Ryzhyk, I.M.} (1980) {\em Table of Integrals, Series and Products}, Acad. Press,
New York.

\item  \textsc{\sc Rainville E.D.} (1960) {\em Special Functions}, Chelsea, New York.

\item  \textsc{\sc Shattuck M., Waldhauser T.}(2010) 'Proofs of some binomial identities using the method of last squares' {\em  Fibonacci Quart.} {\bf 48} 290-297.

\item  \textsc{\sc Sloane N.J.A.} (1964) {\em The On Line Encyclopedia of Integers Sequences }, OEIS.

\item  \textsc{\sc Vervaat W.} (1979) `On a stochastic difference equation and a representation of non-negative infinitely divisible random variables '{\em Adv. Appl. Prob.} {\bf 11} 750-783.

\item  \textsc{\sc WOLFRAM} (2014) {\em Hypergeometric Functions pFq },\\ http://functions.wolfram.com/HypergeometricFunctions/Hypergeometric3F2/13/01/01/01/.

\item  \textsc{\sc WOLFRAM} (2014) {\em On Line Root Finder},\\
http://www.wolframalpha.com/widgets/
\end{enumerate}
\newpage

\section{Appendix: Algorithm for the coefficients}
 \begin{center}

\begin{tabular}{|l|l|l|l|l|l|l| }
\hline

&n-5&n-4&n-3&n-2&n-1&$n=6$ \\
\hline
n&- & -&- &- &- & 1 \\
\hline
n-1&- &- & -&- & 1& $n(n-1)/2=15$\\
\hline
n-2&- &- &- &1&10& $(n-2) table_{n-2,n-1}+table_{n-3,n-1}=65$\\
\hline
n-3&- & -&1&6&25& $(n-3)table_{n-3,n-1}+table_{n-4,n-1}=90$\\
\hline
n-4&- &1&3&7&15& $(n-4)table_{n-4,n-1}+table_{n-5,n-1}=31$\\
\hline
n-5&1&1&1&1&1&  1\\
\hline
* &$\Sigma_5 +$ &$\Sigma_4 +$ &$\Sigma_3 +$ &$\Sigma_2 +$ &$\Sigma_1 +$ &$ \Sigma_0$ \\
\hline
\end{tabular}
\end{center}
Where the $\Sigma$s come from the Viete's formulae for polynoms:\\
$\Sigma_0 = 1$ $ \Sigma_1 = \alpha +\beta +...+\nu $ with n terms
$ \Sigma_2 = \alpha \beta +... $ with $ C_n^2 $ terms
$ \Sigma_k = \alpha \beta... \kappa +... $ with $ C_n^k $ terms
$ \Sigma_n = \alpha \beta ... \nu $ 
Note that this table shows the triangle of Stirling numbers of 2nd kind (see Sloane [13]).

\section{Appendix: Remainder (Cardan) ready for use}
From the non linear system we get the following third order equation, since the corresponding elementary symmetrical functions appear:
$$ \xi^3 -A \xi^2 + (B-A) \xi - C = 0 $$
where:
$$ A  = \alpha+\beta+\gamma+\delta +3 $$
$$ B = \alpha\beta+\alpha\gamma+\alpha\delta+\beta\gamma+\beta\delta+\gamma\delta+
3(\alpha+\beta+\gamma+\delta)+6$$
$$ C = \alpha\beta\gamma+\alpha\gamma\delta+\alpha\beta\delta+\beta\gamma\delta+
\alpha\beta+\alpha\gamma+\alpha\delta+\beta\gamma+\beta\delta+
\gamma\delta+\alpha+\beta+\gamma+\delta+1 $$
According to the classical solution of the third order equation, let:
$$ p = B-A -\frac{A^2}{3}$$
$$ q= -\frac{2 A^3}{27} - C + \frac{A(B-A)}{3}$$
the roots are:
$$ \xi_k = \frac{A}{3} +
j^k \sqrt[3]{-\frac{q}{2} + \sqrt{\frac{q^2}{4}+\frac{p^3}{27}}}
+j^{|k-3|} \sqrt[3]{-\frac{q}{2} - \sqrt{\frac{q^2}{4}+\frac{p^3}{27}}} , k=0,1,2$$
where $j$ is the cubic root of unit.
For instance, if $\alpha=\beta=\gamma=\delta=1$ then $\xi_1 = 3,\xi_2 = 2 + i ,\xi_3 = 2 - i$ Note that the real root $1$ is lacking .

\section{Appendix: Remainder (Ferrari)ready for use}
From the non linear system we get the following 4th order equation, since the corresponding elementary symmetrical functions appear:
$$ \xi^4 -A \xi^3 + (B-3A) \xi^2 - (C-B+2A) \xi + D = 0 $$
where:
$$ A  = \alpha+\beta+\gamma+\delta + \eta +4 $$
$$ B =\alpha\beta+\alpha\gamma+\alpha\delta+\beta\gamma+\beta\delta+\gamma\delta
+\alpha\eta+\beta\eta+\gamma\eta+\delta\eta
+6(\alpha+\beta+\gamma+\delta+\eta)+ 18 $$
\begin{eqnarray*}
C=\alpha\beta\gamma+\alpha\gamma\delta+\alpha\beta\delta+\beta\gamma\delta
+\alpha\beta\eta+\alpha\gamma\eta+\alpha\delta\eta+\beta\gamma\eta+\beta\delta\eta+\gamma\delta\eta\\
+3(\alpha\beta+\alpha\gamma+\alpha\delta+\beta\gamma+\beta\delta+\alpha\eta+\beta\eta+\gamma\eta
+\gamma\delta+\delta\eta)+7(\alpha+\beta+\gamma+\delta+\eta)+14
\end{eqnarray*}
\begin{eqnarray*}
D=\alpha\beta\gamma\eta+\alpha\gamma\delta\eta+\alpha\beta\delta\eta+\beta\gamma\delta\eta
+\alpha\beta\gamma\delta\\
\alpha\beta\gamma+\alpha\gamma\delta+\alpha\beta\delta+\beta\gamma\delta
+\alpha\beta\eta+\alpha\gamma\eta+\alpha\delta\eta+\beta\gamma\eta+\beta\delta\eta+\gamma\delta\eta\\
+(\alpha\beta+\alpha\gamma+\alpha\delta+\beta\gamma+\beta\delta+\alpha\eta+\beta\eta+\gamma\eta
+\gamma\delta+\delta\eta)+(\alpha+\beta+\gamma+\delta+\eta)+1 
\end{eqnarray*}\\
According to the classical solution of the 4th order equation, let:
$$ p = B - 3 A - \frac{3 A^2 }{8} $$
$$ q = \frac{A^3 }{8} +\frac{ A(B-3A)}{2} + C-B + 2 A $$
$$ r = - \frac{3 A^4 }{256} + \frac{(B-3A) A^2 }{16} + \frac{ (C-B+2A)A }{4} + D $$
$$ p_1 = -r -\frac{p^2 }{12}$$
$$ q_1 = - \frac{p^3 }{108} + \frac{4 r p - q^2 }{8} - \frac{p r }{6} $$
the roots of the third order equation are:
$$ \zeta_k = \frac{p}{6} +
j^k \sqrt[3]{-\frac{q_1}{2} + \sqrt{\frac{q_1^2}{4}+\frac{p_1^3}{27}}}
+j^{|k-3|} \sqrt[3]{-\frac{q_1}{2} - \sqrt{\frac{q_1^2}{4}+\frac{p_1^3}{27}}} , k=0,1,2$$
from which $ \zeta_r $ is chosen among the real root(s), and then the 4 roots are:
$$ \xi_k =  \frac{ A }{4} + \frac{1 }{2} \sqrt{2 \zeta_r - p } +(-1)^k \sqrt{\frac{\zeta_r }{2}- \frac{p}{4} +(-1)^{[\frac{k}{2}]} \sqrt{\zeta^2_r - r }} , k=0,1,2,3 $$ For instance,
 if $\alpha=\beta=\gamma=\delta=\eta=1$ then $\xi_1 = 2+ \frac{1}{4}((1-\sqrt{5})+i \sqrt{2(5+\sqrt{5})}),\xi_2 = 2 + \frac{1}{4}((1-\sqrt{5})-i \sqrt{2(5+\sqrt{5})}) ,\xi_3 = 2 + \frac{1}{4}((1+\sqrt{5})+i \sqrt{2(5-\sqrt{5})}), \xi_4 = 2 +\frac{1}{4}((1+\sqrt{5})-i \sqrt{2(5-\sqrt{5})})$ .
 Note that the quasi cyclotomic polynomial has the real root $1$ lacking.
\newpage
\section{ Table of the coefficients of the polynomial for $\alpha=\beta=... =1$}
\begin{center}
\begin{tabular}{|l|l|l|l|l|l|l|l|l|l|l|l|}
\hline
0&1&-& -&- &- &- &-&-&-&-&-\\
\hline
1&1&-3& -&- &- &- &-&-&-&-&-\\
\hline
2&1& -5& 7& -&- &- &- &-&- &-&- \\
\hline
3&1& -7& 17& -15&- &- &- &-&-&-&- \\
\hline
4&1& -9& 31& -49& 31&- &- &- &-&-&-\\
\hline
5&1& -11& 49& -111& 129& -63&- &- &- &-&- \\
\hline
6&1& -13& 71& -209& 351& -321& 127& - &- &-&- \\
\hline
7&1& -15& 97& -351& 769& -1023& 769& -255&-&-&-\\
\hline
8&1& -17& 127& -545& 1471& -2561& 2815& -1793& 511& -&-\\
\hline
9&1& -19& 161& -799& 2561& -5503& 7937& -7423& 4097& -1023&-\\
\hline
10&1&-21 & 199  &-1121 & 4159 &-10625 & 18943  &-23297 & 18943  &-9217 & 2047\\
\hline
n&$x^n$&$x^{(n-1)}$&$x^{(n-2)}$&$x^{(n-3)}$&$x^{(n-4)}$&$x^{(n-5)}$&$x^{(n-6)}$&$x^{(n-7)}$&$x^{(n-8)}$&$x^{(n-9)}$&$x^{(n-10)}$\\
\hline
\end{tabular}
\end{center}

\newpage

\begin{figure}
\centerline{\includegraphics[width=15cm]{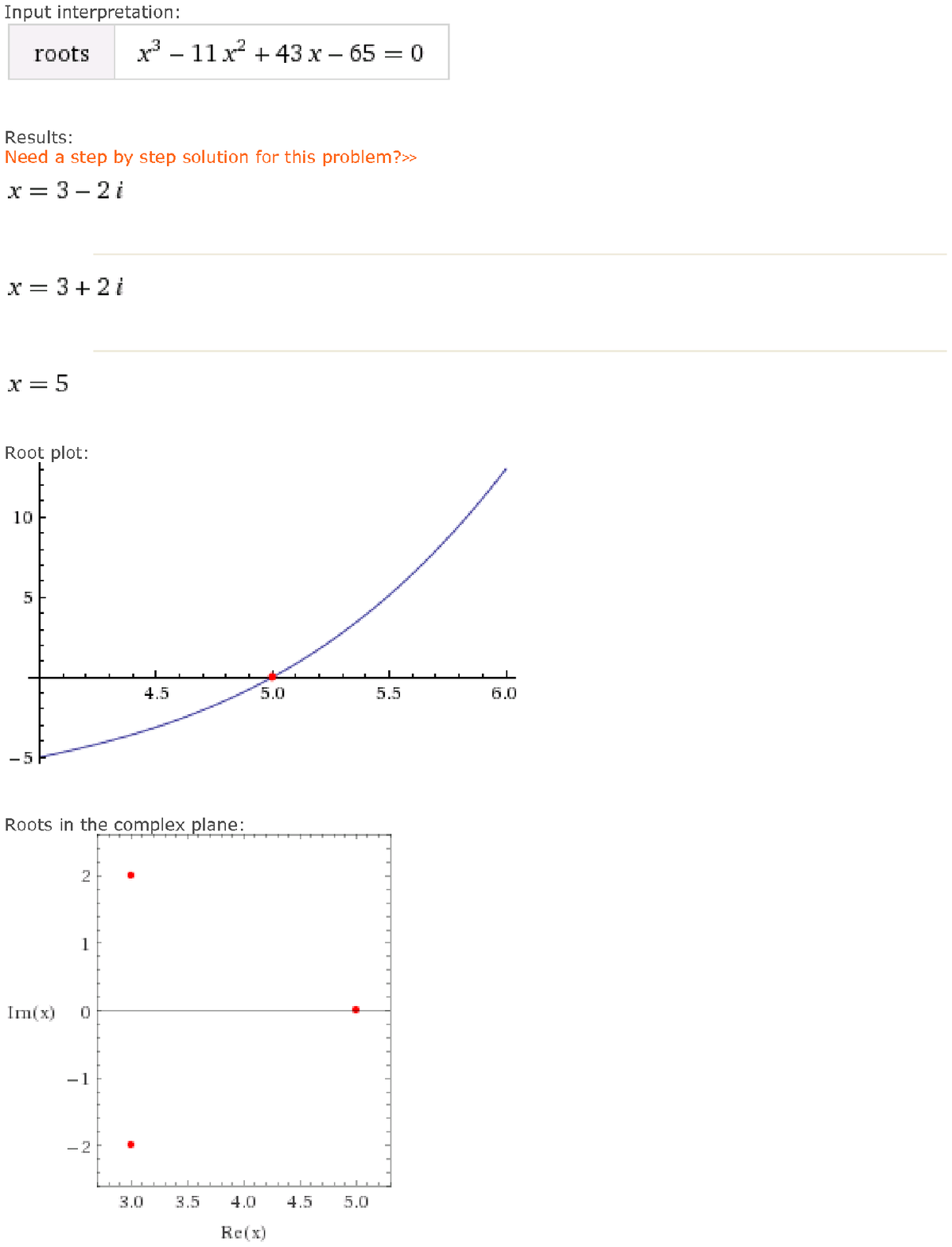}}
\caption{roots for$\alpha= 2$}=
\end{figure} 
\end{document}